\documentclass[10pt,twocolumn,twoside]{IEEEtran}

\IEEEoverridecommandlockouts

\usepackage[cmex10]{amsmath}
\usepackage{fixltx2e}
\usepackage{url}

\usepackage[T1]{fontenc}
\usepackage[latin9]{inputenc}
\usepackage{amstext}
\usepackage{amssymb}
\usepackage{color}

\usepackage{accents}
\usepackage[sort,noadjust]{cite}
\DeclareMathOperator*{\maximize}{\text{maximize}}
\DeclareMathOperator*{\minimize}{\text{minimize}}
\DeclareMathOperator*{\subjectto}{\text{subject to}}

\newtheorem{problem}{Problem}

\newtheorem{definition}{Definition}
\newtheorem{remark}{Remark}

\newtheorem{lemma}{Lemma}
\newtheorem{corollary}{Corollary}
\newtheorem{proposition}{Proposition}
\newtheorem{theorem}{Theorem}
\newtheorem{assumption}{Assumption}
\newcommand{\ubar}[1]{\underaccent{\bar}{#1}}
\newcommand{\norm}[1]{\lVert #1 \rVert}

\newcommand{\abs}[1]{\lvert #1 \rvert}

\DeclareSymbolFont{bbold}{U}{bbold}{m}{n}
\DeclareSymbolFontAlphabet{\mathbbold}{bbold}
\newcommand{\onev}{\mathbbold{1}}

\usepackage{mathtools}

\usepackage[font=small]{caption}
\usepackage[font=small]{subcaption}

\begin{document}

%
\title{Optimal Design of Switched Networks of Positive Linear Systems via Geometric Programming}
%

\author{Masaki~Ogura,~\IEEEmembership{Member,~IEEE,}
        and~Victor~M.~Preciado,~\IEEEmembership{Member,~IEEE}
\thanks{The authors are with the Department
of Electrical and Systems Engineering, University of Pennsylvania, Philadelphia,
PA, 19104 USA. e-mail: \{ogura, preciado\}@seas.upenn.edu.
This work is supported in part by the NSF awards CNS-1302222 and IIS-1447470.
}}
\maketitle

\begin{abstract}
In this paper, we propose an optimization framework to design a network of
positive linear systems whose structure switches according to a Markov process.
The optimization framework herein proposed allows the network designer to
optimize the coupling elements of a directed network, as well as the dynamics of
the nodes in order to maximize the stabilization rate of the network and/or the
disturbance rejection against an exogenous input. The cost of implementing a
particular network is modeled using posynomial cost functions, which allow for a
wide variety of modeling options. In this context, we show that the cost-optimal
network design can be efficiently found using geometric programming in
polynomial time. We illustrate our results with a practical problem in network
epidemiology, namely, the cost-optimal stabilization of the spread of a disease
over a time-varying contact network.
\end{abstract}


%
\IEEEpeerreviewmaketitle

\section{Introduction}

\PARstart{T}{}he intricate structure of many biological, social, and economic
networks emerges as the result of local interactions between agents aiming to
optimize their utilities. The emerging networked system must satisfy both
structural and functional requirements, even in the presence of time-varying
interactions. An important set of functional requirements is concerned with the
behavior of dynamic processes taking place in the network. For example, most
biological networks emerge as the result of an evolutive process that forces the
network to be stable and robust to external perturbations.

{Among networks of dynamical systems, those
consisting of {positive systems} (i.e., the state variables are nonnegative
quantities provided the initial state and the inputs are
nonnegative~\cite{Farina2000}) are of particular importance.} Positive systems arise naturally while modeling
systems in which the variables of interest are inherently nonnegative, such as
concentration of chemical species~\cite{Leenheer2007}, information rates in
communication networks~\cite{Shorten2006}, sizes of infected populations in
epidemiology~\cite{Preciado2014,Ogura2014i}, and many other compartmental
models~\cite{Benvenuti2002}. Due to its practical relevance, many
control-theoretical tools have been adapted to the particular case of positive
linear systems, such as the characterization of stability via diagonal Lyapunov
functions~\cite{Barker1978}, the bounded real lemma~\cite{Tanaka2011}, and
integral linear constraints for robust stability~\cite{Briat2012c}, to mention a
few (we point the reader to~\cite{Farina2000}, for a thorough exposition on
positive systems).

{Most of the methods mentioned above, however, focus on state-feedback control
of positive linear systems. In contrast, there are many situations in which
state-feedback is not a feasible option. A particular example is the
stabilization of a viral spreading process in a complex contact
network~\cite{Preciado2013,Preciado2014}, where it is not feasible to obtain
reliable measurements about the state of the spread in (almost) real time. In
this and similar situations, the control problem is better posed as the problem
of tuning the dynamics of the nodes and the coupling elements of the network in
order to guarantee a stable dynamics. Furthermore, existing control methods
mostly focus on the analysis of time-invariant systems; thus, they do not
provide effective tools for designing time-switching topologies.}

The aim of this paper is to propose a tractable optimization framework to design
networked {positive linear} systems in the presence of time-switching
topologies. We model the time variation of the network structure using Markov
processes, where the modes of the Markov process correspond to different network
structures and the transition rates indicate the probability of switching
between topologies. In this context, we consider the problem of designing both
the dynamics of the nodes in the network, as well as the structure of the
elements coupling them in order to optimize the dynamic performance of the
switching network. We model the cost of implementing a particular network
dynamics using posynomial cost functions, which allow for a wide range of
modeling option \cite{Boyd2007}. We then propose an efficient optimization
framework, based on geometric programming {\cite{Boyd2007}}, to find the
cost-optimal network design that maximizes the stabilization rate and/or the
disturbance rejection against exogenous signals. To achieve this objective, we
develop new theoretical characterizations of the stabilization rate and the
disturbance attenuation of positive Markov jump linear systems which are
specially amenable in the context of geometric programming.

The paper is organized as follows. In Section~\ref{sec:math}, we introduce
elements of graph theory, Markov jump linear systems, and geometric programming
used in our derivations. In Section~\ref{sec:PrbSetting}, we model the dynamics
of randomly switching network using Markov jump linear systems and rigorously
state the design problems under consideration. {Then, we propose geometric
programs to efficiently find the cost-optimal network design for stabilization
and disturbance attenuation, in Sections~\ref{sec:StblOptim}
and~\ref{sec:disturbance}, respectively.} We illustrate our results with a
relevant epidemiological problem, namely, the stabilization of a viral spread in
a switching contact network, in Section~\ref{sec:numerics}.

\section{Mathematical Preliminaries}\label{sec:math}

In this section, we introduce the notation and some results needed in our
derivations. We denote by $\mathbb{R}$ and $\mathbb{R}_+$ the sets of real and
nonnegative numbers, respectively. The set~$\{1,\ldots,N\}$ is sometimes denoted
by $[N]$. We denote vectors using boldface letters and matrices using capital
letters. When $\mathbf{x}\in \mathbb{R}^n$ is nonnegative (positive) entrywise,
we write $\mathbf{x}\geq 0$ ($\mathbf{x}>0$, respectively). The $1$-norm on
$\mathbb{R}^n$ is denoted by $\norm{\mathbf{x}}_1 =\sum_{i=1}^n \abs{x_i}$. The
$n$-vector and $n\times m$ matrix whose entries are all one are denoted by
$\onev_n$ and~$\onev_{n\times m}$, respectively. We denote the identity matrix
by $I$. The $(i, j)$\nobreakdash-entry of a matrix~$A$ is denoted by $A_{i j}$.
The $n\times m$ matrix whose $(i, j)$-entry equals $a_{ij}$ is denoted by
$[a_{ij}]_{i\in[n], j\in [m]}$, or simply by $[a_{ij}]_{i, j}$ when $n$ and $m$
are clear from the context. A square matrix is Hurwitz stable if the real parts
of its eigenvalues are all negative. A square matrix is Metzler if its
off-diagonal entries are nonnegative. The Kronecker product~\cite{Brewer1978} of
the matrices $A$ and $B$ is denoted by $A\otimes B$.  Following
\cite{Asavathiratham2001}, we define the \emph{generalized Kronecker product} of
an $n_1\times n_2$ matrix~$A=[a_{ij}]_{i, j}$ and the set of $m_1\times m_2$
matrices~$\{B_{ij}\}_{i\in[n_1],j\in[n_2]}$ as the following $(n_1 m_1)\times
(n_2 m_2)$ block matrix
\begin{equation*}
A\otimes \{B_{ij}\}_{i, j}
=
\begin{bmatrix}
a_{11}B_{11} & \cdots & a_{1 n_2}B_{1 n_2}\\
\vdots & \ddots & \vdots\\
a_{n_1 1}B_{n_1 1} &  \cdots & a_{n_1 n_2}B_{n_1 n_2}
\end{bmatrix}. 
\end{equation*}
The direct sum of the square matrices $A_1$, $\dotsc$, $A_n$, denoted by
$\bigoplus_{i=1}^n A_i$, is defined as the block-diagonal matrix containing the
matrices $A_1,\dotsc,A_n$ as its diagonal blocks.

\subsection{Basic Graph Theory}

A \emph{weighted}, \emph{directed} graph (also called digraph) is defined as the
triad $\mathcal{G}=\left(\mathcal{V},\mathcal{E},\mathcal{W}\right)$, where
$\mathcal{V}=\left\{ v_{1},\dots,v_{n}\right\} $ is a set of~$n$ nodes,
$\mathcal{E}\subseteq\mathcal{V}\times\mathcal{V}$ is a set of ordered pairs of
nodes called directed edges, and the function
$\mathcal{W}:\mathcal{E}\rightarrow(0, \infty)$ assigns \textit{positive} real
weights to the edges in $\mathcal{E}$. By convention, we say that
$\left(v_{i},v_{j}\right)$ is an edge from $v_{j}$ pointing towards $v_{i}$. The
\emph{adjacency matrix} of a weighted and directed graph $\mathcal{G}$ is an
$n\times n$ matrix $[a_{ij}]_{i, j}$ defined entry-wise as
$a_{ij}=\mathcal{W}((v_{i},v_{j}))$ if edge $(v_{i},v_{j})\in\mathcal{E}$, and
$a_{ij}=0$ otherwise. Thus, the adjacency matrix of a graph is always
nonnegative. Conversely, given a  nonnegative matrix $A$, we can associate to it
a directed graph whose adjacency matrix is $A$.

\subsection{Positive Markov Jump Linear Systems} 

Consider the linear time-invariant system
\begin{equation}\label{eq:LTIsys}
\dot{\mathbf x}= A\mathbf x  + B \mathbf w, 
\quad 
\mathbf z = C\mathbf x  + D \mathbf w, 
\end{equation}
where $\mathbf x(t) \in \mathbb{R}^n$, $\mathbf w(t) \in \mathbb{R}^s$, $\mathbf
z(t) \in \mathbb{R}^r$, and $A$, $B$, $C$, and $D$ are real matrices of
appropriate dimensions. The {system \eqref{eq:LTIsys}} is usually denoted
by the quadruple $(A, B, C, D)$. We say that the {system
\eqref{eq:LTIsys}} is {\it positive} if $\mathbf x(t) \geq 0$ and
$\mathbf z(t)\geq 0$ for all $t\geq 0$ provided that $\mathbf x(0) = \mathbf x_0
\geq 0$ and $\mathbf w(t) \geq 0$ for every $t\geq 0$. It is well
known~\cite{Farina2000} that the {system \eqref{eq:LTIsys}} is positive
if and only if $A$ is Metzler and $B$, $C$, and $D$ are nonnegative.

Let $\sigma = \{\sigma(t)\}_{t\geq 0}$ be a time-homogeneous Markov process with
the state space $\{1, \ldots, M\}$. Let $\Pi =[\pi_{ij}]_{i, j \in [M]}$ denote
the infinitesimal generator of $\sigma$, i.e., assume that the transition
probability of $\sigma$ is given by $\Pr(\sigma(t+h) = j \mid \sigma(t) = i) =
\pi_{ij}h + o(h)$ ($i\neq j$) and $\Pr(\sigma(t+h) = i \mid \sigma(t) = i) =
1+\pi_{ii}h + o(h)$ for all $i, j \in [M]$ and $t, h \geq 0$. We will assume
that the Markov process is irreducible. A Markov jump linear
system~\cite{Costa2013} is described as the following stochastic differential
equations:
\begin{equation}\label{eq:def:Sigma}
\Sigma: \begin{cases}
\dot{\mathbf x}(t) = A_{\sigma(t)}\mathbf x(t) + B_{\sigma(t)}\mathbf w(t),
\\
\mathbf z(t) = C_{\sigma(t)}\mathbf x(t) + D_{\sigma(t)}\mathbf w(t). 
\end{cases}
\end{equation}
We assume that the initial states $\mathbf x(0) = \mathbf x_0$ and $\sigma(0) =
\sigma_0$ are constants. We say that the Markov jump linear system $\Sigma$ is
{\it positive} if the linear time-invariant systems $(A_i, B_i, C_i, D_i)$ are
positive for all $i\in [M]$. We say that  $\Sigma$ is {\it internally mean
stable} (or simply {\it mean stable}) if there exist $\alpha>0$ and $\lambda>0$
such that, if $\mathbf{w}(t) = 0$ for every $t\geq 0$, then $E[\norm{\mathbf
x(t)}_1] \leq \alpha e^{-\lambda t} \norm{\mathbf x_0}_1$ for all $\mathbf x_0$,
$\sigma_0$, and $t\geq 0$, where $E$ denotes the expectation operator. The
supremum of $\lambda$ such that this inequality holds for all $\mathbf x_0$,
$\sigma_0$, and $t$ is called the {\it exponential decay rate} of a mean stable
$\Sigma$. The following proposition describes necessary and sufficient
conditions for mean-stability of positive Markov jump linear systems.

\begin{proposition}[\cite{Bolzern2014,Ogura2013f}]\label{prop:MJLS:stbl}
If the Markov jump linear system~$\Sigma$ is positive, then the following
conditions are equivalent:
\begin{enumerate}
\item $\Sigma$ is mean stable. 

\item The matrix $A = \Pi^\top \otimes I + \bigoplus_{i=1}^M A_i$
is Hurwitz stable.

\item There exist positive vectors $\mathbf{v}_1, \dotsc, \mathbf{v}_M \in
\mathbb{R}^n$ such that $\mathbf{v}_i^\top A_i + \sum_{j=1}^{\raisebox{1.25ex}{}{M}} \pi_{ij}
\mathbf{v}_j^\top < 0$ for every $i\in [M]$.
\end{enumerate}
\end{proposition}

We let $\mathcal{L}_1(\mathbb{R}_+, \mathbb{R}^{n}_+)$ denote the space of
Lebesgue integrable functions on $\mathbb{R}_+$ taking values in
$\mathbb{R}^n_+$. For any $f\in \mathcal{L}_1(\mathbb{R}_+, \mathbb{R}^n_+)$, we
define $\norm{f}_{\mathcal{L}_1} = \int_0^\infty \norm{f(t)}_1\, dt$. The
definition below extends the concept of $\mathcal{L}_1$-stability (given in,
e.g., {\cite{Briat2012c}}) to Markov jump linear systems:

\begin{definition}\label{defn:L1gain}
We say that $\Sigma$ is \emph{$\mathcal{L}_1$-stable} if there exists $\gamma >
0$ such that, for all $\sigma_0$ and $\mathbf w\in \mathcal{L}_1(\mathbb{R}_+,
\mathbb{R}_+^s)$, we have that $E[\mathbf z]\in \mathcal{L}_1(\mathbb{R}_+,
\mathbb{R}_+^r)$ and $\norm{E[\mathbf z]}_{\mathcal{L}_1} < \gamma \norm{\mathbf
w}_{\mathcal{L}_1}$ when $\mathbf x_0 = 0$. If $\Sigma$ is
$\mathcal{L}_1$-stable then its \emph{$\mathcal{L}_1$-gain}, denoted by
$\norm{\Sigma}_1$, is defined by $\norm{\Sigma}_1 = \sup_{\mathbf w\in
\mathcal{L}_1(\mathbb{R}_+, \mathbb{R}_+^s)} ({\norm{E[\mathbf
z]}_{\mathcal{L}_1}}/{\norm{\mathbf w}_{\mathcal{L}_1}})$.
\end{definition}

For the particular case when $\Sigma$ is time-invariant, the next
characterization of the $\mathcal L_1$-gain is available:
\begin{proposition}[\cite{Ebihara2011,Briat2012c}]\label{prop:LTIL1}
{Assume that the linear time-invariant system \eqref{eq:LTIsys} is positive
and let $\gamma>0$ be arbitrary. The following statements are equivalent:
\begin{enumerate}
\item The system {\eqref{eq:LTIsys}} is stable and its $\mathcal{L}_1$-gain
is less than $\gamma$.

\item $A$ is Hurwitz stable and  $\onev^{\!\top}_{r}(D-CA^{-1}B)
- \gamma \onev_s^\top < 0$.

\item There exists a positive vector $\mathbf{v} \in \mathbb{R}^n$ satisfying
inequalities $\mathbf{v}^{\!\top} A + \onev_r^{\!\top} C < 0$ and
$\mathbf{v}^{\!\top} B +  \onev_r^{\!\top} D < \gamma \onev_s^{\!\top}$.
\end{enumerate}}
\end{proposition}

\subsection{Geometric Programming}

The design framework proposed in this paper depends on a class of optimization
problems called geometric programs~\cite{Boyd2007}. Let $x_1$, $\dotsc$, $x_n$
denote $n$ real positive variables and define the vector variable $\mathbf{x} =
(x_1, \dotsc, x_n)$. We say that a real-valued function $g(\mathbf{x})$ is a
{\it monomial function} ({\it monomial} for short) if there exist $c>0$ and
$a_1, \dotsc, a_n \in \mathbb{R}$ such that $g(\mathbf{x}) = c x_{\mathstrut
1}^{a_{1}} \dotsm x_{\mathstrut n}^{a_n}$. Also we say that a real-valued
function~$f(\mathbf{x})$ is a {\it posynomial function} ({\it posynomial} for
short) if it is a sum of monomial functions of $\mathbf{x}$. For information
about the modeling power of posynomials, we point the reader to \cite{Boyd2007}.
Given {positive integers $p, q$ and} a collection of posynomials $f_0{(\mathbf
x)}$, $\dotsc$, $f_p{(\mathbf x)}$ and monomials $g_1{(\mathbf x)}$, $\dotsc$,
$g_q{(\mathbf x)}$, the optimization problem
\begin{equation}\label{eq:scalarGP}
\begin{aligned}
\minimize_{\mathbf x}\ 
&
f_0({\mathbf x})
\\
\subjectto\ 
&
f_i({\mathbf x})\leq 1,\quad i=1, \dotsc, p, 
\\
&
g_j({\mathbf x}) = 1,\quad j=1, \dotsc, q, 
\end{aligned}
\end{equation}
is called a {\it geometric program} (GP) in standard form. Although geometric
programs are not convex, they can be efficiently converted into a convex
optimization problem and efficiently solved using, for example, interior-point
methods (see~\cite{Boyd2007}, for more details on GP). We quote the following
result from \cite[Section~10.4]{Nemirovskii2004} on the computational complexity
of solving GP:

\begin{proposition}\label{prop:GPcomplexity}
{The geometric program \eqref{eq:scalarGP} can be solved with computational
cost $O((k+m)^{1/2} (mk^2 + k^3 + n^3))$, where $m = p+2q$ and $k$ is the
maximum of the numbers of monomials contained in each of posynomials
$f_0(\mathbf x)$, $\dotsc$,~$f_p(\mathbf x)$.}
\end{proposition}

Geometric programs can also be written in terms of matrices and vectors, as
follows. Let $X_1, \dotsc, X_n$ be matrix-valued positive variables and define
$\mathbf{X} = (X_1, \dotsc, X_n)$. We say that {a real matrix-valued function
$F(\mathbf X)$} is a monomial (posynomial) if each entry of $F{(\mathbf X)}$ is
a monomial (respectively, posynomial)  in the entries of the matrix-valued
variables $X_1$, $\dotsc$, $X_n$. Then, given a scalar-valued
posynomial~$f_0{(\mathbf X)}$, matrix-valued posynomials $F_1{(\mathbf X)}$,
$\dotsc$, $F_p{(\mathbf X)}$, and matrix-valued monomials $G_1{(\mathbf X)}$,
$\dotsc$, $G_q{(\mathbf X)}$, we call the optimization problem
\begin{equation*}
\begin{aligned}
\minimize_\mathbf{X}\ 
&
f_0(\mathbf{X})
\\
\subjectto\ 
&
F_i(\mathbf{X})\leq \onev,\quad i=1, \dotsc, p, 
\\
&
G_j(\mathbf{X}) = \onev,\quad j=1, \dotsc, q, 
\end{aligned}
\end{equation*}
a {\it matrix geometric program}, where the right-hand side of the constraints
are all-ones matrices of appropriate dimensions. Notice that one can easily
reduce a matrix geometric program to a standard geometric program by dealing
with the matrix-valued constraints entry-wise.

The next lemma will be used in Sections~\ref{sec:StblOptim} and
\ref{sec:disturbance} to prove our main results.

\begin{lemma}\label{lem:PosyConstraint}
Let $X_0$, $X_1$, $\dotsc$, $X_n$ be matrix-valued positive variables. Assume
that $X_0$ is $\mathbb{R}^{p\times q}$-valued, and let $F(\mathbf{X})$ be an
$\mathbb{R}^{p\times q}$-valued posynomial in the variable $\mathbf{X} = (X_0,
X_1, \dotsc, X_n)$. Then, there exists a matrix-valued posynomial $\tilde
F(\mathbf{X})$ such that $F(\mathbf{X} )\leq X_0$ if and only if $\tilde
F(\mathbf{X} ) \leq \onev_{p\times q}$ for every $\mathbf X$.
\end{lemma}

\begin{IEEEproof}
Define the function $\tilde F$ of $\mathbf{X}$ component-wise as $\tilde
F(\mathbf{X} )_{k, \ell} = (X_{{0}})_{k, \ell}^{-1}F(\mathbf{X} )_{k, \ell}$ for
all $k\in [p]$ and $\ell \in [q]$. Then, $\tilde F(\mathbf{X})$ is a
matrix-valued posynomial. Then, since $X_{{0}}>0$, the constraint $F(\mathbf{X}
)\leq X_{{0}}$ is equivalent to $\tilde F(\mathbf{X} ) \leq \onev_{p\times q}$.
\end{IEEEproof}

\section{Switched Network of Positive Linear Systems} \label{sec:PrbSetting}

In this section, we describe the dynamics of the network of positive linear
systems under consideration and state relevant design problems.

\subsection{Network Dynamical Model} 

Consider a collection of $N$ linear time-invariant subsystems
\begin{equation}\label{eq:subsystems}
\dot{\mathbf{x}}_k = F_k\mathbf{x}_k +G_k\mathbf{u}_k, \quad
\mathbf{y}_k = H_k\mathbf{x}_k,
\end{equation}
where $k \in \mathcal V = \{ 1, \dotsc, N \}$, $\mathbf
x_k(t)\in\mathbb{R}^{n_k}$, $\mathbf y_k(t)\in\mathbb{R}^{p_k}$, and
$\mathbf{u}_k(t) \in \mathbb{R}^{q_k}$. Assume that these subsystems are
linearly coupled through the edges of a {\it time-varying, weighted, and
directed graph} {having a time-dependent adjacency matrix~$K(t) =
[K_{k\ell}(t)]_{k,\ell \in\mathcal V}$}. The coupling between subsystems is
modeled by the following set of inputs:
\begin{equation}\label{eq:coupling}
\mathbf{u}_k(t) = \sum_{\ell \in \mathcal V}K_{k\ell}(t)\Gamma_{\!(k,\ell)}\mathbf{y}_\ell(t), 
\end{equation}
where $k \in \mathcal V$ and $\Gamma_{\!(k,\ell)} \in \mathbb{R}^{q_k\times
p_\ell}$ is the so-called {\it inner coupling matrix}, which indicates how the
output of the $\ell$\nobreakdash-th subsystem influences the input of the
$k$\nobreakdash-th subsystem. The inner-coupling matrix can be interpreted as a
matrix-valued weight associated to each directed edge of the graph.

In the rest of the paper, we consider dynamic networks satisfying the following positivity assumption:

\begin{assumption} 
System~\eqref{eq:subsystems} is positive and the matrix~$\Gamma_{\!(k,\ell)}$ is
nonnegative for all $k, \ell \in \mathcal V$.
\end{assumption}

The above assumption is satisfied for many networked dynamics where the state
variables are nonnegative. For example, the dynamics of many models of disease
spreading in networks~\cite{VanMieghem2009a} can be described as a coupled
network of positive systems, where the state variables represent probabilities
of infection. We will present the details of this model in
Section~\ref{sec:numerics}. Other examples of practical relevance are chemical
networks~\cite{Leenheer2007}, transportation networks~\cite{Rantzer2015}, and
compartmental models \cite{Benvenuti2002}.

In this paper, we consider networks of positive linear systems in which the
network structure switches according to a Markov process, as indicated below:

\begin{assumption}
There exist $N\times N$ nonnegative matrices $K_1$, $\dotsc$, $K_M$ and a
time-homogeneous Markov process \mbox{$\sigma = \{\sigma(t)\}_{t\geq 0}$} taking
its values in $[M]$ such that \mbox{$K(t) = K_{\sigma(t)}$} for every
$t\geq 0$.
\end{assumption}

Under the above assumptions, the dynamics of the network of subsystems in
\eqref{eq:subsystems} coupled through the law \eqref{eq:coupling} forms a
positive Markov jump linear system, as described below. Let us define the
`stacked' vectors $\cramped{\mathbf x = (\mathbf x_1^\top, \dotsc, \mathbf
x_N^\top)^\top}$, \mbox{$\cramped{\mathbf u = (\mathbf u_1^\top, \dotsc, \mathbf
u_N^\top)^\top}$}, and $\cramped{\mathbf y = (\mathbf y_1^\top, \dotsc, \mathbf
y_N^\top)^\top}$. Then, the dynamics in (\ref{eq:subsystems}) can be written by
the differential equations $\dot {\mathbf x } =  \cramped{(\bigoplus_{k=1}^N
F_k)} \mathbf x + \cramped{(\bigoplus_{k=1}^N G_k)} \mathbf u$ and $\mathbf y  =
\cramped{(\bigoplus_{k=1}^N H_k)} \mathbf x$. Similarly, \eqref{eq:coupling} can
be written using the generalized Kronecker product (defined in Section
\ref{sec:math}) as $\mathbf u (t)= ({K(t)}\otimes
\{\Gamma_{\!(k,\ell)}\}_{k,\ell}) \mathbf y (t)$. Thus, the global network
dynamics, denoted by $\mathcal N$, can be written as the following positive
Markov jump linear system:
\begin{equation*} 
\mathcal N: \dot {\mathbf{x}}(t) = A_{\sigma(t)}\mathbf{x}(t),
\end{equation*}
where the matrices $A_i$ ($i=1, \dotsc, M$) are given by
\begin{equation} \label{eq:BigAi}
A_i = \bigoplus_{k=1}^N F_k + \biggl(\bigoplus_{k=1}^N
G_k\biggr)\biggl(K_i\otimes \{\Gamma_{\!(k,\ell)}\}_{k,\ell}\biggr)
\biggl(\bigoplus_{k=1}^N H_k\biggr).
\end{equation}

\subsection{Network Design: Cost and Constraints} 

In this paper, we propose a novel methodology to simultaneously design the
dynamics of each subsystem (characterized by the matrices $F_k$, $G_k$, and
$H_k$ for all $k\in \mathcal V$) and the weights of the edges coupling them
(characterized by the inner-coupling matrices $\Gamma_{\!(k,\ell)}$ for
$k,\ell\in\mathcal V$). Our design objective is to achieve a prescribed
performance criterion (described below) while satisfying certain cost
requirements. In our problem formulation, the time-variant graph structure
changes according to a Markov process $\sigma$ that we assume to be an
\emph{exogenous} signal. In other words, we assume the Markov process ruling the
network switching is out of our control. We refer the readers
to~\cite{Ogura2015d} where, for a model of spreading processes, the authors
study an optimal network design problem over networks whose structure
endogenously changes by reacting to the dynamic state of agents. On the other
hand, we assume that we can modify the dynamics of the subsystems and the matrix
weights of each edge to achieve our design objective.

We specifically assume that, for all $k$ and $\ell$ in $\mathcal V$, there is a
cost function associated to its coupling matrix~$\Gamma_{\!(k,\ell)}$. We denote
this real-valued edge-cost function by $\phi_{(k,\ell)}$. This cost function can
represent, for example, the cost of building an interconnection between two
subsystems. We remark that, in fact, we do not need to design $\Gamma_{\!(k,
\ell)}$ for all $k$ and $\ell$ by the following reason. For each $i\in [M]$, let
$\mathcal G_i = (\mathcal V_i, \mathcal E_i, \mathcal W_i)$ denote the weighted
directed graph having the adjacency matrix $K_i$, and consider the
union~$\mathcal E = \bigcup_{i\in [M]} \mathcal E_i$ of all the possible
directed edges. Then we see that, if $(k, \ell)\notin \mathcal E$, then the
Markov jump linear system~$\mathcal N$ under our consideration is independent of
the value of $\Gamma_{\!(k, \ell)}$. Therefore, throughout the paper, we focus
on designing the coupling matrices $\Gamma_{\!e}$ only for $e\in\mathcal E$.
Similarly, our framework allows us to associate a cost to each subsystem in the
network. For each $k\in \mathcal{V}$, we denote the cost of implementing the
$k$-th subsystem by~$f_k(F_k, G_k, H_k)$. Thus, the total cost of realizing a
particular network dynamics is given by
\begin{equation}
R = \sum_{k\in \mathcal V} f_k(F_k,G_k,H_k) + \,\sum_{\mathclap{e\in \mathcal
E}}\, \phi_{e}(\Gamma_{\!e}). \label{eq:TotalCost}
\end{equation}

In practice, not all realizations of the subsystem $(F_k,G_k,H_k)$ and the
coupling matrix $\Gamma_{\!e}$ are feasible. We account for feasibility
constraints on the design of the $k$-th subsystem via the following set of
inequalities and equalities:
\begin{equation}\label{eq:NodeFeasible}
g_{k,p}^{\mathcal{V}}\left(F_{k},G_{k},H_{k}\right) \leq 1, 
\quad  
h_{k,q}^{\mathcal{V}}\left(F_{k},G_{k},H_{k}\right) =1,
\end{equation}
where $p\in[s_{k}]$, $q\in[t_{k}]$, and $g^{\mathcal{V}}_{k, p}$ and
$h^{\mathcal{V}}_{k, q}$ are real functions for each $k\in\mathcal{V}$.
Similarly, for each edge $e \in\mathcal{E}$, we account for design constraints
on the coupling matrix $\Gamma_{\!e}$ via the restrictions:
\begin{equation}\label{eq:EdgeFeasible}
g_{e,p}^{\mathcal{E}}\left(\Gamma_{\!e}\right) \leq 1,
\quad
h_{e,q}^{\mathcal{E}}\left(\Gamma_{\!e}\right) =1, 
\end{equation}
where $p\in[s_{e}]$, $q\in[t_{e}]$, and $g^{\mathcal{E}}_{e,p}$ and
$h^{\mathcal{E}}_{e,q}$ are real functions for each $e\in\mathcal{E}$. 

\subsection{Network Design: Problem Statements}\label{subset:Problems}

We are now in conditions to formulate the network design problems under
consideration. Our design problems can be classified according to two different
criteria. The first criterion is concerned with the dynamic performance. In this
paper, we limit our attention to two performance indexes: ({\it I}\hspace{.3mm})
{\it stabilization rate} (considered in Section \ref{sec:StblOptim}), and ({\it
II}\hspace{.3mm}) {\it disturbance attenuation} (considered in Section
\ref{sec:disturbance}). The second criterion is concerned with the design cost.
According to this criterion, we have two types of design problems: ({\it A})
{\it performance-constrained} problems and ({\it B}) {\it budget-constrained}
problems. In a budget-constrained problem, the designer is given a fixed budget
and she has to find the network design to maximize a dynamic performance index
(either the stabilization rate or the disturbance attenuation). In a
performance-constrained problem, the designer is required to design a network
that achieves a given performance index while minimizing the total cost of the
design. Combinations of these two criteria described above result in four
possible problem formulations, represented in Table \ref{Table:DesignProblems}.

\begin{table}
\centering
\protect\caption{Design problems under consideration.}
\begin{tabular}{ccc}
\hline 
 & \strut \small Performance-Constr. & \small Budget-Constr. \strut \tabularnewline
\hline 
\small
\strut Stabilization Rate & \small Problem~I-A & \small Problem~I-B \strut
\tabularnewline
\hline 
\strut \small Disturbance Atten. & \small Problem~II-A & \small Problem~II-B \strut
\tabularnewline 
\hline 
\end{tabular}
\label{Table:DesignProblems}
\end{table}

We describe Problems I-A and I-B in more rigorous terms in what follows (Problems II-A and II-B will be formulated in Section~\ref{sec:disturbance}):

\renewcommand{\theproblem}{I-A}
\begin{problem}[Performance-Constrained Stabilization]
Given a desired decay rate ${\lambda} > 0$, design the nodal dynamics
$\{F_k,G_k,H_k\}_{k \in \mathcal V}$ and the coupling matrices
$\{\Gamma_{\!e}\}_{e\in\mathcal E}$ such that the global network dynamics
$\mathcal N$ achieves an exponential decay rate\footnote{We remark that the
exponential decay rate of $\mathcal N$ is well-defined, since $\mathcal N$ is a
Markov jump linear system.} greater than $\lambda$  at a minimum implementation
cost~$R$, defined in \eqref{eq:TotalCost}, while satisfying the feasibility
constraints \eqref{eq:NodeFeasible} and \eqref{eq:EdgeFeasible}.
\end{problem}

\renewcommand{\theproblem}{I-B}
\begin{problem}[Budget-Constrained Stabilization]
Given an available budget $\bar{R} > 0$, design the nodal dynamics
$\{F_k,G_k,H_k\}_{k \in \mathcal V}$ and the coupling matrices
$\{\Gamma_{\!e}\}_{e\in\mathcal E}$ in order to maximize the exponential decay
rate $\lambda$ of $\mathcal N$ while satisfying the budget constraint $R\leq
\bar R$, and the feasibility constraints \eqref{eq:NodeFeasible} and
\eqref{eq:EdgeFeasible}.
\end{problem}

\begin{remark} 
The optimal solutions of Problems~I-A and~I-B are inversely related, as
explained below. Let $R^\star(\lambda)$ ($\lambda^\star(\bar R)$) be the optimal
solution $R$ of Problem~I-A ($\lambda$ of Problem~I\nobreakdash-B,
respectively). For simplicity in our presentation, let us assume that both
$R^\star$ and $\lambda^\star$ are strictly increasing functions and have
nonempty domains. Also, suppose that the compositions~$\lambda^\star \circ
R^\star$ and $R^\star \circ \lambda^\star$ are well-defined. Then, from the
definitions of $\lambda^\star$ and $R^\star$, we have
$\lambda^\star(R^\star(\lambda)) \geq \lambda$ and $R^\star(\lambda^\star(\bar
R)) \leq \bar R$ for all possible values of $\lambda$ and $\bar R$. From these
inequalities, we obtain $\lambda^\star(R^\star(\lambda^\star(\bar R))) \geq
\lambda^\star(\bar R)$ and $\lambda^\star(R^\star(\lambda^\star(\bar R))) \leq
\lambda^\star(\bar R)$, respectively. We therefore have
$\lambda^\star(R^\star(\lambda^\star(\bar R))) = \lambda^\star(\bar R)$. This
yields that $\lambda^\star \circ R^\star$ equals the identity since
$\lambda^\star$ is strictly increasing. In the same way, we can show that
$R^\star \circ \lambda^\star$ is the identity. Hence, $R^\star$ and
$\lambda^\star$ are the inverses of each other.
\end{remark}

In the next section, we proceed to present an optimization framework to solve
Problems I-A and I-B. In Section \ref{sec:disturbance}, we shall extend our
results to Problems II-A and II-B, which we call the Performance-Constrained and
Budget-Constrained Disturbance Attenuation problems, respectively.

\section{Optimal Design for Network Stabilization}\label{sec:StblOptim}

The aim of this section is to present geometric programs in order to solve both the Performance- and the Budget-Constrained Stabilization
problems, which are one of the main contributions of this
paper. In what follows, we place the following
assumption on the cost and constraint functions:

\begin{assumption}\label{assm:posy}\ 
\begin{enumerate}
\item For each $e\in\mathcal E$, the functions $\phi_{e}$ and
$g_{e,p}^{\mathcal{E}}$ ($p\in [s_e]$) are posynomials, while
$h_{e,q}^{\mathcal{E}}$ ($q \in [t_e]$) is a monomial.

\item For each $k\in \mathcal V$, there exists a real and diagonal matrix
{$\Delta_k \in \mathbb{R}^{n_k\times n_k}$} such that, the functions $\hat
f_k(\hat F_k,G_k,H_k) =  f_k(\hat F_k -\Delta_k,G_k,H_k)$ and $\hat
g_{k,p}^{\mathcal{V}}(\hat F_k,G_k,H_k) = g_{k,p}^{\mathcal{V}} (\hat F_k
-\Delta_k,G_{k},H_{k})$ ($p\in[s_{k}]$) are posynomials, while the function
$\hat h_{k,q}^{\mathcal{V}}(\hat F_{k},G_{k},H_{k}) = h_{k,q}^{\mathcal{V}}(\hat
F_{k}-\Delta_k,G_{k},H_{k})$ ($q\in[t_{k}]$) is a monomial.
\end{enumerate}
\end{assumption}

\begin{remark}
Since the state matrix of a positive system is Metzler, it can contain negative
diagonal elements. On the other hand, the decision variables of a geometric
program must be positive. Therefore, negative diagonal entries cannot
necessarily be directly used as decision variables in our optimization
framework. The above assumption will be used to overcome this limitation and
will allow us to design the diagonal elements by a suitable change of variables.
\end{remark}

In order to transform Problem~I-A into a geometric program, we need to introduce
the following definitions. Let 
\begin{equation}\label{eq:delta>=}
\delta 
= 
\max_{i \in [M]} (-\pi_{ii}) 
+ 
\ \ \ 
\max_{\mathclap{k \in \mathcal V,\, i\in [n_k]}}\ \ \ {(\Delta_k)_{ii}}.
\end{equation}
Then, for every $i\in [M]$, define the nonnegative matrix \mbox{$P_i = (\pi_{ii}
+ \delta )I - \bigoplus_{k=1}^N \Delta_k$}. Also, for $\mathbb{R}^{n_k\times
n_k}$-valued positive variables $\hat F_k$ ($k \in \mathcal V$), we define
\begin{equation*}
\hat{A}_i = \bigoplus_{k=1}^N \hat F_k
+ 
\biggl(\bigoplus_{k=1}^N G_k\biggr)
\bigl(K_i\otimes \{\Gamma_{\!(k,\ell)}\}_{k,\ell}\bigr) 
\biggl(\bigoplus_{k=1}^N H_k\biggr). 
\end{equation*}
The next theorem shows how to efficiently solve the Performance-Constrained
Stabilization problem via geometric programming:

\begin{theorem}\label{thm:rate-const-stbl}
The network design that solves Problem~I-A is defined by the set of subsystems
$\{ (F_k^\star, G_k^\star, H_k^\star) \}_{k\in \mathcal V}$ with $F_k^\star =
\hat F_k^\star - \Delta_k$, and the coupling matrices $\{\Gamma_{\!e}^\star
\}_{e\in \mathcal E}$, where the starred matrices are the solutions to the
following matrix geometric program: 
\begin{subequations}
\label{eq:gp:rate-cost-stbl}
\begin{align}
\hspace{-.2cm}\minimize_{\substack{\{\hat F_k, G_k, H_k\}_{k\in\mathcal V},\\ \{\Gamma_{\!e}\}_{e\in\mathcal E}, \{\mathbf{v}_i\}_{i\in[M]}}}
& 
\sum_{{k \in \mathcal V}} \hat f_k(\hat F_k,G_k,H_k)
+
\sum_{{e\in\mathcal E}} \phi_{e}(\Gamma_{\!e})\label{eq:obj}
\\
\text{\upshape subject to }\  \ 
& 
\mathbf{v}_i^{\!\top} (\hat{A}_i +  P_i + {\lambda} I)  + \sum_{{j\neq i}} \pi_{ij}\mathbf{v}_j^{\!\top} < \delta \mathbf{v}_i^{\!\top}\!, \label{eq:stblconst}
\\
& \hat g_{k,p}^{\mathcal{V}}(\hat F_k,G_k,H_k)\leq 1,
\label{eq:const1}
\\
& \hat h_{k,q}^{\mathcal{V}}(\hat F_{k},G_{k},H_{k})=1, 
\label{eq:const2}
\\
&  \eqref{eq:EdgeFeasible}. \label{eq:const3} 
\end{align}
\end{subequations}
\end{theorem}

\begin{remark}\label{rem:GPRemark}
The optimization program \eqref{eq:gp:rate-cost-stbl} is, in fact, a matrix
geometric program. The cost function in \eqref{eq:obj} is a posynomial under
Assumption \ref{assm:posy}. Also, the set of design constraints
\eqref{eq:const1}--\eqref{eq:const3} are valid posynomial inequalities and
monomial equalities. Furthermore, the constraint in \eqref{eq:stblconst} is a
matrix-posynomial constraint by Lemma~\ref{lem:PosyConstraint}. Also the
definition of $\delta$ in \eqref{eq:delta>=} ensures that the matrix $P_i$ is
nonnegative.
\end{remark}

\begin{remark}
Standard GP solvers cannot handle strict inequalities, such as
\eqref{eq:stblconst}. In practice, we can overcome this limitation by including
an arbitrary small number to relax the strict inequality into a non-strict
inequality.
\end{remark}

Before we present the proof of Theorem \ref{thm:rate-const-stbl}, we need to
introduce the following corollary of Proposition~\ref{prop:MJLS:stbl}.

\begin{corollary}\label{cor:}
Assume that the Markov jump linear system~$\Sigma$ defined
in~\eqref{eq:def:Sigma} is positive. Then, $\Sigma$ is mean stable with an
exponential decay rate greater than $\lambda>0$, if and only if there exist
positive vectors $\mathbf{v}_1, \dotsc, \mathbf{v}_M \in \mathbb{R}^n$ such that
\mbox{$\mathbf{v}_i^\top A_i + \sum_{j=1}^M \pi_{ij} \mathbf{v}_j^\top + \lambda
\mathbf{v}_i^\top < 0$} for every $i \in [M]$.
\end{corollary}

\begin{IEEEproof} 
Let us assume that $\Sigma$ is mean stable and its exponential decay rate is
greater than $\lambda>0$. Then, the Markov jump linear system $\dot{\mathbf x}
(t) = (A_{\sigma(t)}+\lambda I) \mathbf x(t)$ is mean stable. Therefore, by
Proposition~\ref{prop:MJLS:stbl}, we can find positive vectors $\mathbf{v}_i \in
\mathbb{R}^n$ such that $\mathbf{v}_i^\top (A_i + \lambda I) + \sum_{j=1}^M
\pi_{ij} \mathbf{v}_j^\top < 0$. This proves the necessity of the condition in
the corollary. The sufficiency part can be proved in a similar way.
\end{IEEEproof}

Let us prove Theorem~\ref{thm:rate-const-stbl}.

\begin{IEEEproof}[Proof of Theorem~\ref{thm:rate-const-stbl}]
By Corollary~\ref{cor:}, the Performance-Constrained Stabilization problem is
equivalent to the following optimization problem:
\begin{subequations}
\label{eq:gp:rate-cost-stbl:pf}
\begin{align}
\hspace{-.27cm}
\minimize_{\substack{\{F_k, G_k, H_k\}_{k\in \mathcal V},\\
\{\Gamma_{\!e}\}_{e\in \mathcal E},
\{\mathbf{v}_i\}_{i\in[M]}}}
& 
\sum_{k\in \mathcal V} f_k(F_k,G_k,H_k) 
+ 
\,\sum_{{e\in \mathcal E}}\, \phi_{e}(\Gamma_{\!e}) \label{eq:minimize_R}
\\
\text{\upshape subject to}\  \ \ 
& \mathbf{v}_i^{\!\top} A_i + \sum_{j=1}^M \pi_{ij}\mathbf{v}_j^{\!\top}
+ \lambda \mathbf{v}_i^{\!\top}< 0,\, \mathbf{v}_i > 0,
\label{eq:pf:stblconst}
\\
& \text{\eqref{eq:NodeFeasible} and \eqref{eq:EdgeFeasible}}. 
\end{align}
\end{subequations}
Notice that $A_i$, defined in \eqref{eq:BigAi}, can have negative diagonal
entries, since $F_k$ is a Metzler matrix. If this is the case, the constraint in
\eqref{eq:pf:stblconst} cannot be written as a posynomial inequality. To
overcome this issue, for each $k \in \mathcal V$, we use the following
transformation
\begin{equation}\label{eq:F_k=...}
\hat F_k = F_k + \Delta_k, 
\end{equation}
where $\Delta_k$ is given in Assumption~\ref{assm:posy}. In fact, if the
triple~$(F_k, G_k, H_k)$ is a feasible solution of
\eqref{eq:gp:rate-cost-stbl:pf}, then $\hat f_k(\hat F_k, G_k, H_k) = f_{k}(F_k,
G_k, H_k)$ is well-defined, which shows that $\hat F_k$ is positive because
$\hat f_k$ is a posynomial.

Then, we show that \eqref{eq:pf:stblconst} is equivalent to
\eqref{eq:stblconst}. Noting that the transformation \eqref{eq:F_k=...} yields
$\hat A_i = A_i + {\bigoplus_{k=1}^N \Delta_k}$\vspace{.1\baselineskip}, we can
rewrite the constraint \eqref{eq:pf:stblconst} as $\mathbf{v}_i^\top \hat A_i
-\mathbf{v}_i^\top {\bigoplus_{k=1}^N \Delta_k} + \pi_{ii} \mathbf{v}_i^\top+
\sum_{j\neq i} \pi_{ij} \mathbf{v}_j^\top + \lambda \mathbf{v}_i^\top <
0$\vspace{.1\baselineskip}. Adding $\delta \mathbf{v}_i^\top$ to both sides of
the above inequality, we obtain \eqref{eq:stblconst}. Also, the equivalence
between the objective functions~\eqref{eq:minimize_R} and \eqref{eq:obj} is
obvious from their definitions. In the same way, we can observe that the
constraints~\eqref{eq:NodeFeasible} and \eqref{eq:EdgeFeasible} are equivalent
to the constraints \eqref{eq:const1}\nobreakdash--\eqref{eq:const3}. Therefore,
we conclude that the optimization problems \eqref{eq:gp:rate-cost-stbl:pf} and
\eqref{eq:gp:rate-cost-stbl} are equivalent. {Notice that the constraint
$\mathbf{v}_i > 0$ is omitted in \eqref{eq:gp:rate-cost-stbl} because
\eqref{eq:gp:rate-cost-stbl} is stated as a geometric program and, therefore,
all of its variables are assumed to be positive.} Finally, as mentioned in
Remark~\ref{rem:GPRemark}, the optimization problem in
\eqref{eq:gp:rate-cost-stbl} is indeed a matrix geometric program.
\end{IEEEproof}

Theorem \ref{thm:rate-const-stbl} allows us to find the cost-optimal network
design to stabilize the system at a desired decay rate, assuming the feasible
set defined by \eqref{eq:stblconst}--\eqref{eq:const3} is not empty. Similarly,
the following theorem introduces a geometric program to solve the
Budget-Constrained Stabilization problem:

\begin{theorem}\label{thm:bud-const-stbl}
The network design that solves Problem~I-B is described by the set of subsystems
$\{ (F_k^\star, G_k^\star, H_k^\star) \}_{k\in \mathcal V}$ with $F_k^\star =
\hat F_k^\star - \Delta_k$, and the coupling matrices $ \{\Gamma_{\!e}^\star
\}_{e\in \mathcal E}$, where the starred matrices are the solutions to the
following matrix geometric program:
\begin{equation} \label{eq:gp:bud-const-stbl}
\begin{aligned}
\hspace{-.25cm}
\minimize_{\substack{\lambda,\{\hat F_k, G_k, H_k\}_{k\in \mathcal V},\\ \{\Gamma_{\!e}\}_{e\in \mathcal E}, \{\mathbf{v}_i\}_{i\in[M]}}}
& 
1/\lambda\\
\text{\upshape subject to}\ \ \ \ 
& 
\text{\eqref{eq:stblconst}--\eqref{eq:const3}, }
\\
& 
\sum_{k\in\mathcal V} \hat f_k(\hat F_k,G_k,H_k)
+ 
\sum_{e\in\mathcal E} \phi_{e}(\Gamma_{\!e})\leq\bar R.
\end{aligned}
\end{equation}
\end{theorem}

\begin{IEEEproof}
By Corollary~\ref{cor:}, we can formulate the Budget-Constrained
Stabilization problem as the following optimization problem:
\begin{align}
\maximize_{\substack{\lambda,\{F_k, G_k, H_k\}_{k\in \mathcal V},\\ \{\Gamma_{\!e}\}_{e\in \mathcal E},\{\mathbf{v}_i\}_{ i\in[M]}}}\  
& 
\lambda
\label{eq:optim:budsta}
\\
\text{\upshape subject to}\ \ \ \ \ 
& 
\text{\eqref{eq:NodeFeasible}, \eqref{eq:EdgeFeasible}, \eqref{eq:pf:stblconst}, 
$\mathbf v_i > 0$, $\lambda > 0$, and $R \leq \bar R$.}\notag
\end{align}
Notice that maximizing $\lambda$ is equivalent to minimizing $1/\lambda$, since
$\lambda >0$. The rest of the proof is similar to the proof of
Theorem~\ref{thm:rate-const-stbl}, and we limit ourselves to remark the main
steps. First, we use the transformation in~\eqref{eq:F_k=...} to show that the
optimization problem~\eqref{eq:optim:budsta} is equivalent to
\eqref{eq:gp:bud-const-stbl}. We can show that the optimization
problem~\eqref{eq:gp:bud-const-stbl} is a matrix geometric program by following
the explanation in Remark~\ref{rem:GPRemark}.
\end{IEEEproof}

\begin{remark}
Using a discrete-time counterpart of Proposition~\ref{prop:MJLS:stbl}, one could
easily establish discrete-time analogues of Theorems~\ref{thm:rate-const-stbl}
and \ref{thm:bud-const-stbl} in a straightforward manner (see, e.g., the
stability characterizations of discrete-time positive Markov jump linear systems
in \cite[Theorem~1]{Zhu2014a} and \cite[Theorem~3.4]{Ogura2013f}).
\end{remark}

\section{Optimal Design for Disturbance Attenuation}\label{sec:disturbance}

In Section \ref{sec:StblOptim}, we have introduced an optimization framework,
based on geometric programming, to solve both the Performance- and
Budget-Constrained Stabilization problems described in
Subsection~\ref{subset:Problems}. In this section, we extend our analysis to
design networks from the point of view of disturbance attenuation. We consider
the following collection of linear time-invariant subsystems:
\begin{equation}\label{eq:subsys_w/_disturbance}
\begin{aligned}
\dot{\mathbf{x}}_k& = F_k\mathbf{x}_k+ G_{k,1} \mathbf w_k + G_{k,2} \mathbf{u}_k,
\\
\mathbf z_k &= H_{k,1} \mathbf x_k + J_{k,11} \mathbf w_k + J_{k,12} \mathbf{u}_k,
\\
\mathbf{y}_k& = H_{k,2}\mathbf{x}_k+ J_{k,21} \mathbf w_k, 
\end{aligned}
\end{equation}
where $k \in \mathcal V$. The signals $\mathbf w_k(t)\in \mathbb{R}^{s_k}$ and
$\mathbf z_k(t)\in\mathbb{R}^{r_k}$ are, respectively, the disturbance input and
the performance output. We assume that subsystems are linearly interconnected
through the law \eqref{eq:coupling}, and that the following positivity
assumption holds:

\begin{assumption}
For every $k \in \mathcal V$, $F_k$ is Metzler and $G_{k,1}$, $G_{k,2}$,
$H_{k,1}$, $H_{k,2}$, $J_{k,11}$, $J_{k,12}$, and $J_{k,21}$ are nonnegative.
\end{assumption}  

{Using the generalized Kronecker product, we can rewrite the dynamics of
the network $\mathcal N$ of subsystems in \eqref{eq:subsys_w/_disturbance},
coupled according to~\eqref{eq:coupling}, as a Markov jump linear system
\eqref{eq:def:Sigma} with the following coefficient matrices:
\begin{equation*}
\begin{aligned}
A_i &= \bigoplus_{k=1}^N F_k + \biggl(\bigoplus_{k=1}^N G_{k,2}\!\biggr) (K_i\!\otimes \!\{\Gamma_{\!(k,\ell)}\}_{k,\ell}) \biggl(\bigoplus_{k=1}^N H_{k,2}\!\biggr), 
\\
B_i &= \bigoplus_{k=1}^NG_{k,1} + \biggl(\bigoplus_{k=1}^N G_{k,2}\!\biggr) (K_i\!\otimes \!\{\Gamma_{\!(k,\ell)}\}_{k,\ell}) \biggl(\bigoplus_{k=1}^NJ_{k,21}\!\biggr), 
\\
C_i &= \bigoplus_{k=1}^N H_{k,1} + \biggl( \bigoplus_{k=1}^N J_{k,12}\!\biggr) (K_i\!\otimes\! \{\Gamma_{\!(k,\ell)}\}_{k,\ell}) \biggl(\bigoplus_{k=1}^N H_{k,2}\!\biggr), 
\\
D_i &= \bigoplus_{k=1}^N J_{k,11} + \biggl(\bigoplus_{k=1}^N J_{k,12}\!\biggr) (K_i\!\otimes \!\{\Gamma_{\!(k,\ell)}\}_{k,\ell}) \biggl(\bigoplus_{k=1}^N J_{k,21}\!\biggr). 
\end{aligned}
\end{equation*}
Therefore, the $\mathcal{L}_1$-gain {(see Definition~\ref{defn:L1gain})} of
the switched network of linearly coupled subsystems is well defined. We can then
rigorously state Problem~II-A, as follows: \renewcommand{\theproblem}{II-A}
\begin{problem}[Performance-Constrained Disturbance Attenuation]
Given a desired $\mathcal{L}_1$-gain, denoted by $\gamma > 0$, design the nodal
dynamics $\{F_k,G_k,H_k\}_{k\in\mathcal V}$ and the coupling
matrices~$\{\Gamma_{\!e}\}_{e\in\mathcal E}$ such that the global network
dynamics $\mathcal N$ is mean stable and its $\mathcal{L}_1$-gain is less than
$\gamma$, while the implementation cost $R$ defined in \eqref{eq:TotalCost} is
minimized and the feasibility constraints \eqref{eq:NodeFeasible} and
\eqref{eq:EdgeFeasible} are satisfied.
\end{problem}}

The next theorem, which gives the second main contribution of this paper,
shows that this problem can be solved via geometric programming:
\begin{theorem}\label{thm:L1-GainConstrained}
{The network design that solves Problem II-A is defined by the set of subsystems
$\{ (F_k^\star, G_k^\star, H_k^\star) \}_{k\in \mathcal V}$ with $F_k^\star =
\hat F_k^\star - \Delta_k$, and the coupling matrices $ \{\Gamma_{\!e}^\star
\}_{e\in \mathcal E}$, where the starred matrices are the solutions to the
following matrix geometric program:
\begin{align*}
\minimize_{\substack{\lambda,\{\hat F_k, G_k, H_k\}_{k\in\mathcal V},\\ \{\Gamma_{\!e}\}_{e\in\mathcal E}, \{\mathbf{v}_i\}_{i\in[M]}}}\ \notag
& 
\sum_{\mathclap{k\in\mathcal V}} \hat f_k(\hat F_k,G_k,H_k)
+
\sum_{\mathclap{e\in\mathcal E}} \phi_{e}(\Gamma_{\!e})\notag
\\
\text{\upshape subject to}\ \ \ \ \ 
& 
\mathbf{v}_i^{\!\top}(\hat{A}_i + P_i)  + \sum_{j\neq i} \pi_{ij}\mathbf{v}_j^{\!\top}  + \onev_r^{\!\top} C_i < \delta \mathbf{v}_i^{\!\top},\notag
\\
&
\mathbf{v}_i^\top B_i + \onev_r^\top D_i < \gamma \onev_s^\top,\ 
\\
& 
\text{\eqref{eq:const1}--\eqref{eq:const3}}.
\end{align*}}
\end{theorem}

\subsection{Proof}

The proof of Theorem~\ref{thm:L1-GainConstrained} is based on  the
following theorem, which reduces the analysis of the $\mathcal{L}_1$-gain of a
general positive Markov jump linear system to a linear program\footnote{A
preliminary version of the theorem can be found in \cite{Ogura2014g}. Also, we
remark that the theorem is a continuous-time counterpart of
\cite[Theorem~2]{Zhu2014a}.}:

\begin{theorem}\label{thm:L1char}
Assume that the Markov jump linear system~$\Sigma$ defined
in~\eqref{eq:def:Sigma} is positive. For any $\gamma>0$, $\Sigma$  is internally
mean stable and $\norm{\Sigma}_1 < \gamma$, if and only if there exist positive
vectors~$\mathbf{v}_1, \dotsc, \mathbf{v}_M \in \mathbb{R}^{n}$ such that
\begin{equation}\label{eq:LP}
\mathbf{v}_i^{\!\top} A_i 
+ 
\sum_{j=1}^M \pi_{ij}\mathbf{v}_j^{\!\top} 
+ 
\onev_r^{\!\top}  C_i 
< 0,
\ 
\mathbf{v}_i^{\!\top} B_i + \onev_r^{\!\top} D_i < \gamma \onev_s^{\!\top}, 
\end{equation}
for every $i\in \{1, \dotsc, M\}$. 
\end{theorem}

For the proof of Theorem~\ref{thm:L1char}, we introduce $\{0, 1\}$-valued
stochastic processes $\xi_1$, $\dotsc$, $\xi_M$, defined as $\xi_i(t) = 1$ if
\mbox{$\sigma(t) = i$} and $\xi_i(t) = 0$ otherwise for every $i$ and $t$, and
define the vector $\xi = (\xi_1,\dotsc,\xi_M)^\top$. Following the steps in the
proof of \cite[Proposition~5.3]{Ogura2013f}, one can prove that:
\begin{equation}\label{eq:invariantized}
\begin{aligned}
\frac{d}{dt}E[\xi \otimes \mathbf x] 
&= 
 A E[\xi \otimes \mathbf x]  +  B (E[\xi]\otimes \mathbf w), 
\\
E[\xi \otimes \mathbf z] 
&= 
 C E[\xi\otimes \mathbf x] +  D (E[\xi]\otimes \mathbf w), 
\end{aligned}
\end{equation}
where $A$ is defined in Proposition~\ref{prop:MJLS:stbl}, $ B=
\bigoplus_{i=1}^M B_i$, $C = \bigoplus_{i=1}^M C_i$, and $ D =
\bigoplus_{i=1}^M D_i$. We can then prove the following useful lemma:

\begin{lemma} 
If $\Sigma$ is positive and internally mean stable, $\mathbf x_0 = 0$, and
$\mathbf w(t)\geq 0$ for every $t\geq 0$, then 
\begin{equation}\label{eq:lem:tf(0)}
\int_0^\infty E[\xi\otimes
\mathbf z] \,dt = (D - C A^{-1} B ) \int_0^\infty (E[\xi]\otimes \mathbf
w)\,dt.
\end{equation}
\end{lemma}

\begin{IEEEproof}
Assume that $\Sigma$ is internally mean stable and set $\mathbf x_0 = 0$. Then,
$A$ is Hurwitz stable by Proposition~\ref{prop:MJLS:stbl}. From
\eqref{eq:invariantized}, it follows that $E[\xi(t)\otimes \mathbf z(t)] =
\int_0^t \! C e^{ A(t-s)} B (E[\xi(s)]\otimes \mathbf w(s))\,ds + D
(E[\xi(t)]\otimes \mathbf w(t))$. Since $A$ is Hurwitz stable and, also, the
functions $C e^{ A(t-s)} B$, $E[\xi\otimes \mathbf z]$, and $E[\xi]\otimes
\mathbf w$ are nonnegative, integrating this equation from $0$ to $\infty$ with
respect to $t$ completes the proof.
\end{IEEEproof}

We also prove the following lemma, which can provide alternative expressions for
$\norm{\mathbf w}_{\mathcal{L}_1}$ and $\norm{E[\mathbf z]}_{\mathcal{L}_1}$:

\begin{lemma}\label{lem:lifts}
Assume that $\Sigma$ is positive. If $\mathbf x_0 = 0$ and $\mathbf w(t) \geq 0$
for every $t\geq 0$, then
\begin{align} \label{eq:||w||}
\int_0^\infty \! \norm{\mathbf w(t)}_1\,dt
&=
\onev_{Ms}^\top 
\int_0^\infty \! (E[\xi] \otimes \mathbf w)\,dt, 
\\
\int_0^\infty \! \norm{E[\mathbf z(t)]}_1\,dt
&=
\onev_{Mr}^\top
\int_0^\infty \!E[\xi \otimes \mathbf z]\,dt. \label{eq:||E[z]||}
\end{align}
\end{lemma}

\begin{IEEEproof} 
We only give the proof of the second equation~\eqref{eq:||E[z]||}, since the
proof of the first one is identical. From the assumptions in the statement of
the lemma, we have that $\mathbf z(t) \geq 0$ for every $t\geq 0$ with
probability one. Therefore, by the linearity of expectations, the identity
$\norm{\mathbf v}_1 = \onev_{n}^\top \mathbf v$ that holds for a general
$\mathbf v\in\mathbb{R}^n_+$, and the fact that $\onev^{\!\top}_M\xi = 1$, we
can show that $\onev_{Mr}^{\!\top}E[\xi(t)\otimes \mathbf z(t)] =
E[(\onev_{M}^{\!\top}\xi(t))(\onev_{r}^{\!\top}\mathbf z(t))] =
\onev_{r}^{\!\top}E[\mathbf z(t)] = \norm{E[\mathbf z(t)]}_1$. Integrating the
both sides of this equation with respect to $t$ from $0$ to $\infty$, we obtain
\eqref{eq:||E[z]||}.
\end{IEEEproof}

We now have the elements needed to prove Theorem~\ref{thm:L1char}. 

\begin{IEEEproof}[Proof of Theorem~\ref{thm:L1char}]
First assume that $\Sigma$ is internally mean stable and $\norm{\Sigma}_1
<\gamma$. Then, $A$ is Hurwitz stable by Proposition~\ref{prop:MJLS:stbl} and,
therefore, invertible. Thus, the vector~$\eta = \onev_{Mr}^\top ( D - C^\top
A^{-1} B) - \gamma  \onev_{Ms}^\top$ is well-defined. Let us first show $\eta<
0$. Take $\epsilon > 0$ such that $\norm{\Sigma}_1 <\gamma-\epsilon$.  Then, for
every initial state $\sigma_0$ and $\mathbf w \in \mathcal L_1(\mathbb{R}_+,
\mathbb{R}_+^s)$, if $\mathbf x_0=0$, then $E[\mathbf z] \in \mathcal
L_1(\mathbb{R}_+, \mathbb{R}_+^r)$ and $\norm{E[\mathbf z]}_{\mathcal{L}_1} <
(\gamma-\epsilon) \norm{\mathbf w}_{\mathcal{L}_1}$. Therefore, by
Lemma~\ref{lem:lifts}, we have that $\onev_{Mr}^\top \int_0^\infty \!
E[\xi\otimes \mathbf z]\,dt < (\gamma-\epsilon) \onev_{Ms}^\top \int_0^\infty
(E[\xi]\otimes \mathbf w) \, dt$. Then, by \eqref{eq:lem:tf(0)}, we obtain
\begin{equation}\label{eq:(eta...<0}
(\eta+ \epsilon \onev_{Ms}^\top)  
\int_0^\infty \! (E[\xi] \otimes \mathbf w)\,dt < 0. 
\end{equation}
Now, let $i \in [M]$ and~$j \in [s]$ be arbitrary. Let $\mathbf e_i$ and
$\mathbf f_j$ be the $i$-th and $j$-th standard unit vectors in $\mathbb{R}^M$
and~$\mathbb{R}^s$, respectively. Let $\sigma_0 = i$ and, for every $\tau > 0$,
define $\mathbf w_\tau = \tau\chi_{[0, 1/\tau]}\mathbf f_j \in
\mathcal{L}_1(\mathbb{R}_+, \mathbb{R}^s_+)$, where $\chi_{[0, 1/\tau]} \colon
\mathbb{R} \to \{0, 1\}$ denotes the indicator function of the interval $[0,
1/\tau]$. Notice that, by a standard argument in distribution
theory~\cite{Zemanian1965}, the function $\tau \mapsto \tau \chi_{[0, 1/\tau]}$
converges to the Dirac delta function as $\tau \to \infty$ in the space of
distributions. Thus, in the limit of \mbox{$\tau \to \infty$}, we obtain
\mbox{$\int_0^\infty (E[\xi] \otimes \mathbf w_\tau) \,dt \to (E[\xi(0)])
\otimes \mathbf f_j =\mathbf  e_i \otimes \mathbf f_j$}. Therefore,
\eqref{eq:(eta...<0} shows that $(\eta + \epsilon \onev_{Ms}^\top) (\mathbf
e_i\otimes \mathbf f_j) \leq 0$. Since $i$ and~$j$ are arbitrary, it must be the
case that $\eta + \epsilon \onev_{Ms}^\top\leq 0$ and, hence, $\eta < 0$, as we
wanted to show. Now, by Proposition~\ref{prop:LTIL1}, there exists a positive
\mbox{$\mathbf{v} \in \mathbb{R}^{Mn}$} such that
\begin{equation}\label{eq:twoLPs}
\mathbf{v}^{\!\top}  A + \onev_{Mr}^\top  C < 0,
\quad
\mathbf{v}^{\!\top}  B +  \onev_{Mr}^\top  D <  \gamma \onev_{Ms}^\top . 
\end{equation}
Then, from the definition of the matrix $A$ (see
Proposition~\ref{prop:MJLS:stbl}), it is easy to see that the inequalities
\eqref{eq:LP} are satisfied by the positive vectors $\mathbf{v}_1, \dotsc,
\mathbf{v}_M \in \mathbb{R}^n$ given by $\mathbf{v} = (\mathbf{v}_1^\top,
\hdots, \mathbf{v}_M^\top)^\top$. This completes the proof of the necessity part
of the theorem.

On the other hand, assume that there exist positive vectors $\mathbf{v}_1,
\dotsc, \mathbf{v}_M \in \mathbb{R}^n$ such that \eqref{eq:LP} holds. Then we
can see that the positive vector $\mathbf{v} = (\mathbf{v}_1^\top, \dotsc,
\mathbf{v}_M^\top)^\top\in\mathbb{R}^{Mn}$ satisfies \eqref{eq:twoLPs}.
Therefore, by Proposition~\ref{prop:LTIL1}, the linear time-invariant positive
system $(A, B, C, D)$ is stable {and has the $\mathcal L_1$-gain less than
$\gamma$}. Thus, $A$ is Hurwitz stable and, by Proposition~\ref{prop:MJLS:stbl},
the Markov jump linear system $\Sigma$ is internally mean stable. We need to
show that $\norm{\Sigma}_1 < \gamma$. Let $\sigma_0$ and $\mathbf w \in
\mathcal{L}_1(\mathbb{R}_+, \mathbb{R}_+^s)$ be arbitrary and set $\mathbf x_0 =
0$. Let $\mathbf z$ denote the corresponding trajectory of $\Sigma$. From
\eqref{eq:||E[z]||} and \eqref{eq:lem:tf(0)}, it follows that $\int_0^\infty
\norm{E[\mathbf z(t)]}_1\,dt = \onev^{\top}_{Mr}(D-CA^{-1}B) \int_0^\infty
(E[\xi]\otimes \mathbf w)\,dt$. Since the positive linear time-invariant system
$(A, B, C, D)$ has the $\mathcal L_1$-gain less than $\gamma$,
Proposition~\ref{prop:LTIL1} shows that there exists an $\epsilon>0$ satisfying
\mbox{$\onev^{\!\top}_{Mr}(D-CA^{-1}B)< (\gamma-\epsilon) \onev^{\!\top}_{Ms}$}.
Therefore, we obtain $\int_0^\infty \norm{E[\mathbf z(t)]}_1\,dt <
(\gamma-\epsilon) \norm{\mathbf w}_{\mathcal L_1}$, where we used
\eqref{eq:||w||}. This proves $E[\mathbf z] \in \mathcal L_1(\mathbb{R}_+,
\mathbb{R}_+^s)$ and therefore the $\mathcal L_1$-stability of $\Sigma$.
Moreover we have $\norm{\Sigma}_1\leq \gamma-\epsilon < \gamma$. This completes
the proof of the theorem.
\end{IEEEproof}

Now we can give the proof of Theorem~\ref{thm:L1-GainConstrained}:

\begin{IEEEproof}[Proof of Theorem~\ref{thm:L1-GainConstrained}]
From Theorem~\ref{thm:L1char}, the Performance-Constrained Disturbance
Attenuation problem can be stated as the following optimization problem:
\begin{align*}
\minimize_{{\lambda,\{\hat F_k, G_k, H_k\}_{k\in \mathcal V}, \{\Gamma_{\!e}\}_{e\in\mathcal E}, \{\mathbf{v}_i\}_{i\in[M]}}}\ 
& 
R \notag
\\
\text{\upshape subject to}\ \ \ \ \ \ \ \ \ \ \ \ \ \,
& 
\text{\eqref{eq:NodeFeasible}, \eqref{eq:EdgeFeasible}, $\mathbf{v}_i> 0$, and \eqref{eq:LP}.}
\end{align*}
Applying the change of variables in~\eqref{eq:F_k=...}, we obtain the
optimization program stated in the theorem. We can also see that the
optimization problem is, in fact, a matrix geometric program by following the
explanation in Remark \ref{rem:GPRemark}.
\end{IEEEproof}

\begin{remark}
As mentioned in Subsection \ref{subset:Problems}, Problem~II\nobreakdash-B
considers the case of minimizing $\gamma$ (i.e., maximizing the disturbance
attenuation) under a budget constraint. The statement and solution of
Problem~II-B fall straightforward from the cases previously considered and,
hence, details are omitted.
\end{remark}

\section{Numerical Simulations}\label{sec:numerics} 

In this section, we illustrate our network design framework to stabilize the
dynamics of a disease spreading in a time-varying network of individuals. We
consider a popular networked dynamic model from the epidemiological literature,
the networked Susceptible-Infected-Susceptible (SIS)
model~\cite{VanMieghem2009a}. According to this model, the evolution of the
disease in a networked population can be described as:
\begin{equation}\label{eq:N-inter}
\dot{x}_k(t) = -\delta_k x_k(t) + \beta_k \sum_{\ell \in \mathcal V} K_{k\ell}(t)
x_\ell(t) + \epsilon_k w_k(t)
\end{equation}
for $k \in \mathcal V$, where $x_k(t)$ is a scalar variable representing the
probability that node $k$ is infected at time $t$. The parameter~$\delta_k>0$,
called the \emph{recovery rate}, indicates the rate at which node $k$ would be
cured from a potential infection. The parameter $\beta_k>0$, called the
\emph{infection rate}, indicates the rate at which the infection is transmitted
to node $k$ from its infected neighbors. The exogenous signal $\epsilon_k
w_k(t)$, where $\epsilon_k\geq 0$ is a constant and $w_k$ is an
$\mathbb{R}_+$-valued function, is introduced to explain possible transmission
of infection from outside of the network. The entries of the time-varying
adjacency matrix of the contact network are $K_{k\ell}(t)\in\{0,1\}$, for $k,
\ell \in \mathcal V$.

We consider the following epidemiological problem \cite{Preciado2014}: Assume we
have access to vaccines that can be used to reduce the infection rates of
individuals in the network, as well as antidotes that can be used to increase
their recovery rates. Assuming that both vaccines and antidotes have an
associated cost, how would you distribute vaccines and antidotes throughout the
individuals  in the network in order to eradicate an epidemic outbreak at a
given exponential decay rate while minimizing the total cost? We state this
question in rigorous terms below and present an optimal solution using geometric
programming. Let $c_1(\beta_k)$ and $c_2(\delta_k)$ denote the costs of tuning
the infection rate $\beta_k$ and the recovery rate $\delta_k$ of agent $k$,
respectively. We assume that these rates can be tuned within the following
feasible intervals:
\begin{equation}\label{eq:gammabetabounds}
0< \ubar \beta_k \leq \beta_k \leq \bar \beta_k, \ 
0< \ubar \delta_k \leq
\delta_k \leq \bar \delta_k. 
\end{equation}
The problem of finding the optimal allocation of vaccines and antidotes in a
static network was recently solved in \cite{Preciado2014}, under certain
assumptions on the cost functions $c_1$ and $c_2$. In what follows, we solve
this problem for {\it time-varying} contact networks, $K(t) =
K_{\sigma(t)}$, where~$\sigma$ is a time-homogeneous Markov process. In
particular, when there is no disturbance from outside of the network, i.e., when
$\epsilon_k=0$ for every $k\in \mathcal V$, we can formulate this problem as a
particular version of the Performance-Constrained Stabilization problem solved
in Section \ref{sec:StblOptim}, as follows:

\renewcommand{\theproblem}{1}
\begin{problem}\label{prb:epidemiology}
Assume $\epsilon_k=0$ for every $k\in \mathcal V$. Given a desired decay rate
$\lambda > 0$, tune the spreading and recovery rates~$\{\beta_k\}_{k\in \mathcal
V}$ and $\{\delta_k\}_{k\in \mathcal V}$ in the network such that the disease
modeled in~\eqref{eq:N-inter} is eradicated at an exponential decay rate of
$\lambda$ and a minimum cost $R = \sum_{k \in \mathcal V} (c_1(\beta_k) +
c_2(\delta_k))$, while satisfying the box constraints in
\eqref{eq:gammabetabounds}.
\end{problem}

We can transform the above problem into a Performance-Constrained Stabilization
problem, as follows. It is easy to see that the systems \eqref{eq:N-inter} form
a network of positive linear systems $\dot{x}_k = -\delta_k x_k + \beta_k u_k+
\epsilon_k w_k$ with the coupling $u_k(t) = \sum_{\ell\in\mathcal V}
K_{k\ell}(t) \beta_k x_\ell(t)$. We set $f_k(F_{k},G_k,H_k) = c_2(-F_k)$ for
each $k\in \mathcal V$ and $\phi_{(k,\ell)} =c_1/d_k$ for each $(k, \ell)
\in\mathcal E$, where $d_k$ is the indegree of vertex $k$ defined by \mbox{$d_k
= \abs{\{\ell\in\mathcal V: (k, \ell)\in\mathcal E\}}$}. For illustration
purposes, we use the set of cost functions proposed in~\cite{Preciado2014}:
\begin{equation*}
c_1(\beta_k) = \frac{\beta_k^{-1} -{\bar \beta}_k^{-1}}
{{\ubar \beta}_k^{-1} - {\bar \beta}_k^{-1}}, \ 
c_2(\delta_k) =\frac{(1-\delta_k)^{-1} - (1-\ubar \delta_k)^{-1}}
{(1-\bar \delta_k)^{-1}-(1-\ubar \delta_k)^{-1}}. 
\end{equation*}
Notice that $c_1$ is decreasing, $c_2$ is increasing, and the range of $c_1$ and
$c_2$ are both $[0, 1]$. In this particular example, we let {$\Delta_k =
1$} for every $k \in \mathcal V$. Then, it follows that
\begin{equation*}
\hat f_k(\hat F_k,G_k,H_k) 
=
\frac{\hat F_k^{-1} - (1-\ubar \delta_k)^{-1}}{
(1-\bar \delta_k)^{-1}-(1-\ubar \delta_k)^{-1}}.
\end{equation*}
The negative constant term $-(1-\ubar \delta_k)^{-1}$ in the numerator of $\hat
f_k$ can be ignored, since it only changes the value of the total cost $R$ by a
fixed constant (while the optimal allocation is unchanged). For the same reason,
we can ignore the negative constant $-{\bar \beta}_k^{-1}$ in the numerator of
$\phi_{(k,\ell)}$. Notice that, after ignoring these constant terms, the
functions $\hat f_k$ and $\phi_{(k,\ell)}$ are posynomials and, hence,
Assumption~\ref{assm:posy} holds true. Notice also that, since $\beta_k$ does
not depend on $\ell$, we need the additional constraints $\Gamma_{\!(k,1)} =
\cdots = \Gamma_{\!(k,N)}$, which can be implemented as monomial constraints
$\Gamma_{\!(k,1)}/\Gamma_{\!(k,\ell)} = 1$ for $\ell = 2, \dotsc, N$.

Once the cost functions are defined, let us consider a particular model of
network switching to illustrate our results. Consider a set of agents $1,
\dotsc, N$ divided into $h$ disjoint subsets called {\it Households}. Each
Household has exactly one agent called a {\it Worker}. The set of Workers is
further divided into $w$ disjoint subsets called {\it Workplaces}. We assume
that the topology of the contact network switches between three possible graphs.
In the first contact graph $K_1$, a pair of agents are adjacent if they are
in the same household. The second graph $K_2$ is a random graph in which any
pair of Workers are adjacent with probability $p = 0.3$, independently of other
interactions. We also assume that non-Worker agents are adjacent if they belong
to the same Household. The third graph~$K_3$ is a collection of disconnected
$h+w$ complete subgraphs. In particular, we have $w$ complete subgraphs formed
by disjoint sets of Workers sharing the same Workplace, and $h$ complete
subgraphs formed by sets of non-Workers belonging to the same Household.

We consider the case in which the topology of the network switches according to
a Markov process $\sigma$ taking its values in the set $\{1, 2, 3, 2\}$ with the
following infinitesimal generator:
\begin{equation}\label{eq:Pi}
\Pi = \begin{bmatrix}
-1/13  &  1/13  &  0 & 0\\
0 &  -1 & 1 &  0\\
0& 0  & -1/9 & 1/9\\
1 & 0 & 0 & -1
\end{bmatrix}. 
\end{equation}
This choice reflects Workers' simplified schedule: 13 hours of stay at
home (6 PM--7 AM), 1 hours of commute (7 AM--8 AM), 9 hours of
work (8 AM--5 PM), and 1 hour of commute (5 PM--6 PM). 

In this setup, we find the optimal allocation of resources to control an
epidemic outbreak for the following parameters. First, we randomly generate a
set of $n=247$ agents with $h=71$ Workers (as many as Households) and $w=10$
Workplaces. We let $\ubar \beta_k = 0.01$, $\bar \beta_k = 0.05$, $\ubar
\delta_k = -0.5$, $\bar \delta_k = -0.1$, and $\epsilon_k=0$ for every $k\in
\mathcal V$. We then solve the Performance-Constrained Stabilization problem
(Problem~\ref{prb:epidemiology}) with $\lambda = 0.01$ using
Theorem~\ref{thm:rate-const-stbl} and obtain the optimal investment strategy
with total cost of $R = 54.38$.

\begin{figure}[tb]
\centering \includegraphics[width=.97\linewidth]{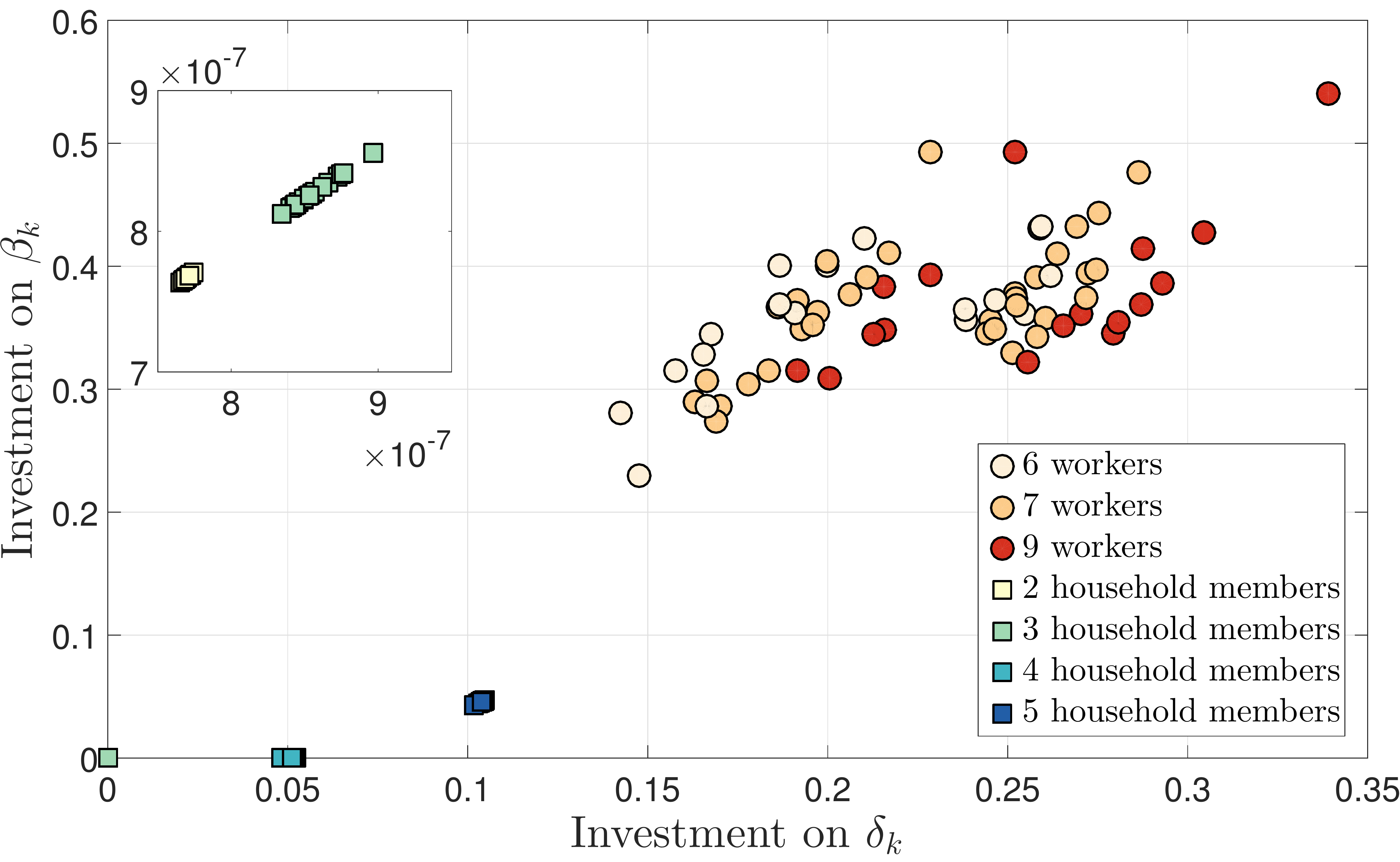}
\caption{Correction versus prevention per node. Circle: Workers. Square:
non-Workers. The color of markers indicate the size of Workplaces/Households
where agents belong to.} \label{fig:1}
\end{figure}

\begin{figure}[tb]
\vspace{.1cm}
\begin{minipage}[b]{\linewidth}
\centering
\includegraphics[width=.97\linewidth]{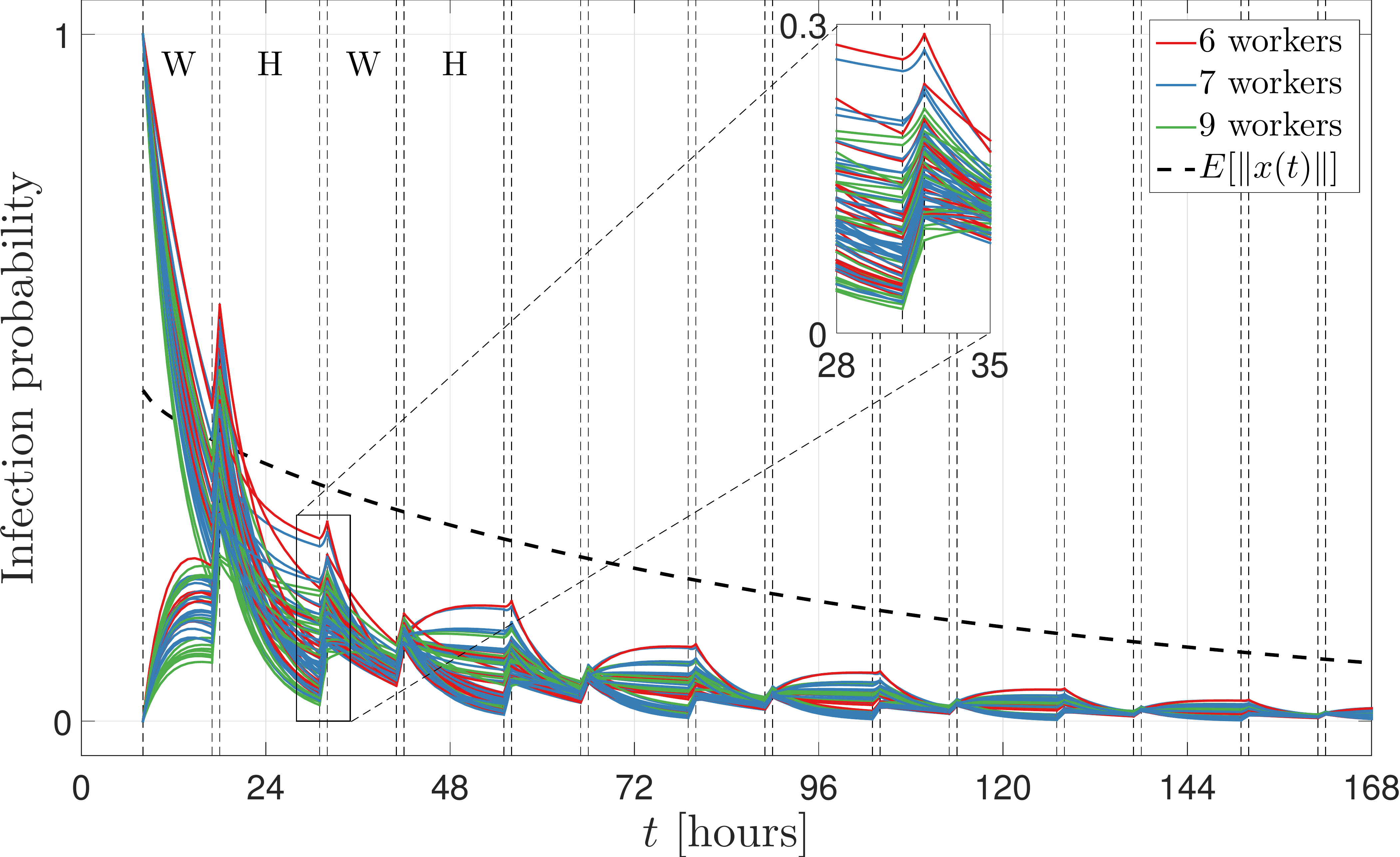}
\subcaption{Infection probabilities of Workers.
}\label{fig:2home}
\end{minipage}
\vspace{.2cm}
\\
\begin{minipage}[b]{\linewidth}
\centering
\includegraphics[width=.97\linewidth]{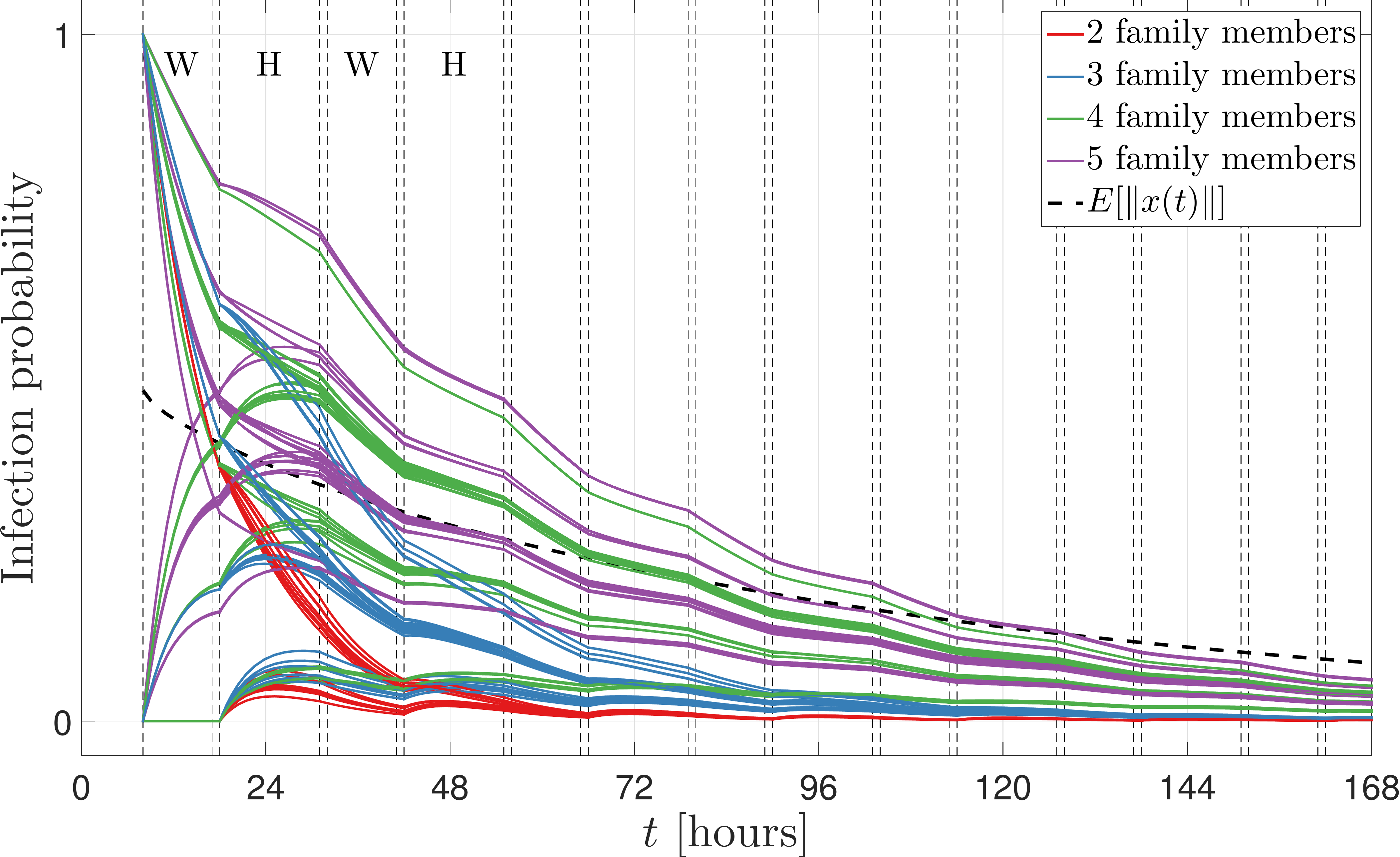}
\subcaption{Infection probabilities of non-Workers. 
}\label{fig:2comm}
\end{minipage}
\caption{Solid: Infection probabilities $x_k(t)$ of Workers (\subref{fig:2home})
and non-Workers (\subref{fig:2comm}). Dotted: $E[\norm{x(t)}_1]$. The colors of
the lines indicate the size of Workplaces/Households where agents belong to. The
labels 'W' and 'H' indicate the first two time-periods when Workers are in their
Workplaces and Households, respectively. Narrow vertical stripes represent the
time-periods when Workers are commuting.}
\end{figure}

Fig.~\ref{fig:1} is a scatter plot showing the relationship between the
investments on $\delta_k$ (corrective actions) and $\beta_k$ (preventive
actions) for each agent. We can see that, in order to maximize the effectiveness
of our budget, we need to invest heavily on Workers, in particular on those
belonging to larger Workplaces. Figs.~\ref{fig:2home} and \ref{fig:2comm} show
the infection probabilities~$x_k$ starting from {$t=t_0 = 8$}, when Workers
start working {(therefore $\sigma(t_0) = 3$)}. We choose the vector of initial
probabilities of infection, {$x_k(t_0)$}, at random from the set $\{0, 1\}^N$.
Each solid curve in the figures represents the evolution in the probability of
infection of an agent over time. The dashed vertical lines indicate the times
when the graph changes. We are also including a dashed curve representing
$E[\norm{x(t)}]$ when $\sigma$ follows the Markov process with the infinitesimal
generator~\eqref{eq:Pi}. 

\begin{figure}[tb]
\centering \includegraphics[width=.97\linewidth]{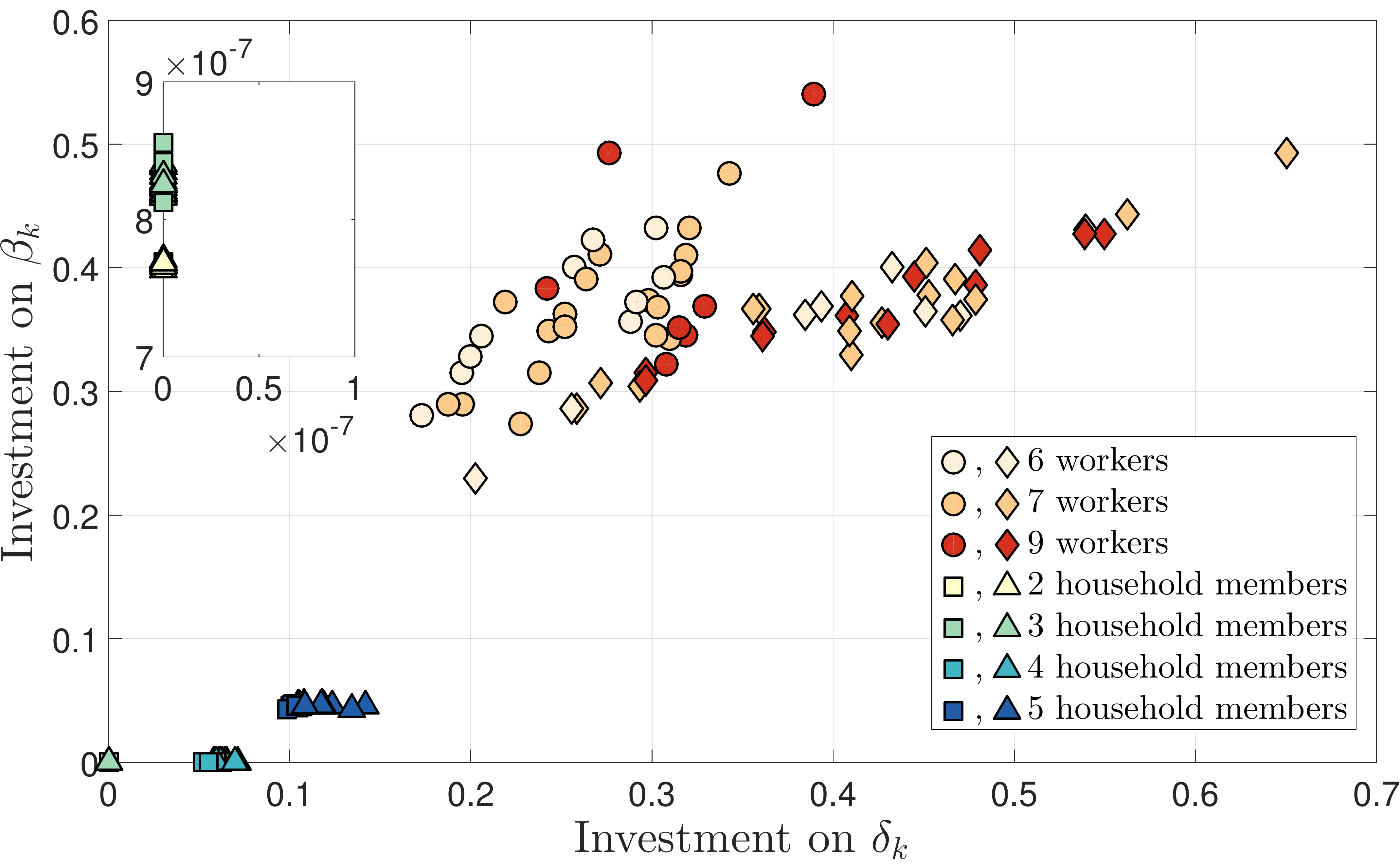}
\caption{Correction versus prevention per node. Circle:  Workers. Diamond:  Workers with disturbance. Square:  non-Workers. Triangle:  non-Workers with disturbance. The color of markers indicate the size of Workplaces/Households
where agents belong to.} \label{fig:3}
\end{figure}

In the second example, we randomly set $\epsilon_k$ to be either $0$ or $1$ for
each $k\in \mathcal V$, and solve the Performance-Constrained Disturbance
Attenuation problem with the constraint $\sup_{w\in\mathcal L_1(\mathbb{R}_+,
\mathbb{R}_+^N)} (\norm{E[\mathbf x]}_{\mathcal L_1}/\norm{\mathbf w}_{\mathcal
L_1}) < \gamma = 40$. We use the same cost functions as in our first example.
Using Theorem~\ref{thm:L1-GainConstrained}, we obtain the optimal values of
$\beta_k$ and $\delta_k$ with a total cost of $67.81$. Fig.~\ref{fig:3} is a
scatter plot showing the relationship between the investments on $\delta_k$ and
$\beta_k$. We can observe different patterns of resource allocation for agents
with and without disturbance; in general, those with disturbance receive more
allocation for corrective resources, while those without disturbance do more for
preventive resources.

Before closing this section, we briefly discuss the computational cost of
solving the geometric programs in this example. For simplicity, we assume that
the number of nonzero entries in each column of $\Pi$ is less than or equal to
$N$. Then, according to the notation in Proposition~\ref{prop:GPcomplexity}, we
can show that $n = N(M+2)$, $k \leq 2N$, and $m \leq N(M+4)$ for
Problem~\ref{prb:epidemiology}. Therefore, by
Proposition~\ref{prop:GPcomplexity}, the computational cost for solving
Problem~\ref{prb:epidemiology} equals $O(N^{7/2}M^{7/2})$. Similarly, we can
show that the same computational cost is required by the Performance-Constrained
Disturbance Attenuation problem.

\section{Conclusions {and Discussion}}

In this paper, we have proposed an optimization framework to design the
subsystems and coupling elements of a time-varying network to satisfy certain
structural and functional requirements. We have assumed there are both
implementation costs and feasibility constraints associated with these network
elements, which we have modeled using posynomial cost functions and
inequalities, respectively. In this context, we have studied several design
problems aiming at finding the cost-optimal network design satisfying certain
budget and performance constraints, in particular, stabilization rate and
disturbance attenuation. We have developed new theoretical tools to cast these
design problems into geometric programs and illustrated our approach by solving
the problem of stabilizing a viral spreading process in a time-switching contact
network.

A possible direction for future research is the extension to
networks of nonpositive systems. The proposed framework is not
directly applicable to this case due to the positivity constraint in geometric
program. However, by applying the framework to the upper-bounding linear
dynamics of vector Lyapunov functions presented in~\cite{Ikeda1980}, we may be
still able to accomplish an efficient, if not optimal, network design for
stabilization and disturbance attenuation.


\begin{IEEEbiography}{Masaki Ogura} received his B.Sc.~degree in Engineering and M.Sc.~degree in Informatics from Kyoto University, Japan, in 2007 and 2009, respectively, and his Ph.D.~degree in Mathematics from Texas Tech University in 2014. He is currently a Postdoctoral Researcher
in the Department of Electrical and Systems Engineering at the University of Pennsylvania. His research interest includes dynamical systems on time-varying 
networks, switched linear systems, and stochastic processes.
\end{IEEEbiography}

\begin{IEEEbiography}{Victor M. Preciado} received the Ph.D. degree in Electrical Engineering
and Computer Science from the Massachusetts Institute of Technology,
Cambridge in 2008.

He is currently the Raj and Neera Singh Assistant Professor of Electrical
and Systems Engineering at the University of Pennsylvania. He is a member of the Networked and Social Systems Engineering (NETS) program and the Warren Center for Network and Data Sciences. His research interests include network science, dynamic systems, control theory, and convex optimization with applications in socio-technical networks, technological infrastructure, and biological systems.
\end{IEEEbiography}







\end{document}